\documentclass{article}

\usepackage{amsfonts}
\usepackage{amsmath}
\usepackage{amssymb}
\input{xy}
\xyoption{all}

\providecommand{\acht}{\scriptstyle}

\providecommand{\rad}{\mathop{\rm rad}\nolimits}
\providecommand{\Hom}{\mathop{\rm Hom}\nolimits}
\providecommand{\Ext}{\mathop{\rm Ext}\nolimits}
\providecommand{\End}{\mathop{\rm End}\nolimits}
\providecommand{\modd}{\mathop{\rm mod}\nolimits}

\providecommand{\soc}{\mathop{\rm soc}\nolimits}
\providecommand{\tp}{\mathop{\rm top}\nolimits}
\providecommand{\im}{\mathop{\rm Im}\nolimits}
\providecommand{\Ker}{\mathop{\rm Ker}\nolimits}

\providecommand{\isomorphe}{\cong}

\providecommand{\CA}{\mathcal{C}_A}
\providecommand{\DA}{\mathcal{D}^b( \textup{mod}\, A  )}
\providecommand{\DC}{\mathcal{D}^b( \textup{mod}\, C  )}

\providecommand{\zd}{\delta}

\providecommand{\zg}{\gamma}

\providecommand{\za}{\alpha}
\providecommand{\zb}{\beta}
\providecommand{\ze}{\epsilon}
\providecommand{\zl}{\lambda}

\newtheorem{theorem}{Theorem}[section] 
\newtheorem{lemma}[theorem]{Lemma}

\newtheorem{definition}[theorem]{Definition}
\newtheorem{remark}[theorem]{Remark}
\newtheorem{example}[theorem]{Example}
\newenvironment{proof}{{\it Proof}.}

\title{Cluster-tilted algebras as trivial extensions}
\author{I. Assem, T. Br\"ustle and R. Schiffler\footnote{The first and the second author gratefully acknowledge partial support
  form the NSERC of Canada. The second author also thanks the universities of
 Sherbrooke and Bishop's for partial support. The third author was
  partially supported by the University of Massachusetts.}}
\date{}
\begin{document}
\maketitle

  \begin{abstract} \noindent Given a  finite
dimensional algebra $C$ (over an algebraically closed field) of global
dimension  at most two, we define its relation-extension algebra to be the trivial extension $C\ltimes
\Ext_C^2(DC,C)$ of $C$  by the $C$-$C$-bimodule $\Ext_C^2(DC,C)$. We
give a construction for the quiver of the relation-extension algebra in case the
quiver of $C$ has no oriented cycles. Our main result says that an
algebra $\tilde C$ is cluster-tilted if and
    only if there exists a tilted algebra $C$ such that $\tilde C$ is
    isomorphic to the relation-extension of $C$.
\end{abstract}

\begin{section}{Introduction} 
Cluster categories were introduced in  \cite{BMRRT}, and, for type
$A_n$ also in \cite{CCS}, as a means for a better understanding of the
cluster algebras of Fomin and Zelevinsky \cite{FZ1,FZ2}. They are
defined as follows: let $A$ be
a hereditary algebra, and $\DA$ be the  derived
category of bounded complexes of finitely generated $A$-modules, then
the cluster category $\CA$ is the orbit category of  $\DA$ 
under the action of the functor $F=\tau^{-1}\,[1]$,
where $\tau$ is the Auslander-Reiten translation in
$\DA $ and $[1]$ is the shift.

In \cite{BMR1},  Buan, Marsh and Reiten defined the cluster-tilted
algebras as follows. Let $A$ be a hereditary algebra, and
$\tilde T$ be a tilting object in $\mathcal{C}_A$, that is,  an object
such that 
$\Ext^1_{\CA}(\tilde T,\tilde T)=0$ and 
the number of isomorphism classes of  indecomposable summands of
$\tilde T$ equals the number of isomorphism classes of simple $A$-modules.
Then the endomorphism algebra  $\End_{\CA}(\tilde T)$ is called
cluster-tilted. 
Since then, these algebras have been the subject of many
investigations, see, for instance, \cite{BMR1,BMR2,BMR3,BR,BRS,CCS,CCS2,KR}.
In several particular cases, it was shown that the quiver of a
cluster-tilted algebra  was obtained from that of a tilted algebra by
replacing relations by arrows, see, for instance
\cite{BR,BRS}. 
Our objective in this paper is to prove this statement in a more
general context (not depending on the representation type).  
This is achieved by looking at cluster-tilted algebras as trivial
extensions of tilted algebras by a bimodule which we explicitely
describe (compare \cite{B}).

For this purpose, we let $C$ be a finite dimensional algebra of global
dimension two (over an algebraically closed field), and consider the
$C$-$C$-bimodule $\Ext_C^2(DC,C)$ with the natural action. The
trivial extension $C\ltimes \Ext_C^2(DC,C)$ is called the
relation-extension algebra of $C$. Our first main result (Theorem
\ref{2.5}) describes the quiver of the relation-extension of $C$ in
the case where the quiver of $C$ has no oriented cycles: we prove that
indeed this quiver is given by replacing each element in a (minimal)
system of relations by an
arrow (going in the opposite direction to the relation). We then prove
the main result of this paper.
\begin{theorem}
  An algebra $\tilde C$ is cluster-tilted if and
    only if there exists a tilted algebra $C$ such that $\tilde C$ is
    the relation-extension of $C$.
\end{theorem} 

We note that several tilted algebras may correspond to the same
cluster-tilted algebra, so this mapping is not bijective. On the other
hand, there clearly exist relation-extension algebras which are not
cluster-tilted.

Combining the above theorem with  Theorem \ref{2.5} we deduce the
construction of the quiver of a cluster-tilted algebra. This allows,
for instance, as done in \cite{BRS}, to relate the list of tame
concealed algebras of Happel and Vossieck \cite{HV} with Seven's list
of minimal infinite cluster quivers \cite{S}.

This paper consists of two sections. The first one describes
relation-extension algebras and their quivers, and the second is
devoted to the cluster-tilted algebras. Moreover, we give several examples.

Th. Br\"ustle wishes to thank Claus Michael Ringel and Idun Reiten for interesting discussions on this problem.
\end{section}
\begin{section}{Relation-extension algebras}
\begin{subsection}{The definition} Throughout this paper, algebras are basic and
  connected finite dimensional algebras over a fixed algebraically
  closed field $k$. For an algebra $C$, we denote by $\modd
  C $ the category of finitely generated right $C$-modules and by
  $\DC$ the derived category of bounded complexes over $\modd C$. The
  functor $D=\Hom_k(-,k)$ is the standard duality between $\modd C$
  and $\modd C^{op}$. For facts about $\modd C$ or $\DC$, we refer to
  \cite{ARS,R,H}.

Let  $C$ be an algebra. We recall that the \emph{trivial extension} of
  $C$ by a $C$-$C$-bimodule $M$ is the algebra $C\ltimes M$ with
  underlying $k$-vector space
\[C\oplus M = \{(c,m) \mid c\in C, m\in M\} \]
and the multiplication defined by
\begin{equation}\nonumber
  (c,m)(c',m')=(cc',cm'+mc')
\end{equation}
for $c,c'\in C$ and $m,m'\in M$. For trivial
extension algebras, we refer to \cite{FGR,AM}.

In this section, we introduce a particular class of trivial extension
algebras which are useful for studying the cluster-tilted algebras.

\begin{definition} Let $C$ be a finite dimensional algebra of global
  dimension at most two, and consider the $C$-$C$-bimodule
  $\Ext_C^2(DC,C)$ (with the natural action). The trivial extension
  \[ C\ltimes \Ext_C^2(DC,C)\]
  is called the \emph{relation-extension} of $C$.
\end{definition}

Clearly, any hereditary algebra is (trivially) the relation-extension
of itself. On the other hand, if $C$ is of global dimension equal to
two (thus not hereditary) there exist two simple $C$-modules $S$ and $S'$
such that $\Ext^2_C(S,S')\ne 0$. Denoting by $I$ the injective
envelope of $S$ and by $P'$ the projective cover of $S'$, the short
exact sequences 
\[\xymatrix{0\ar[r] &S\ar[r] & I \ar[r] &I/S\ar[r] &0}\] 
\[\textup{and} \quad \xymatrix{ 0\ar[r] &\rad P' \ar[r] &P'\ar[r] &S' \ar[r] &0} \]
induce an epimorphism $\Ext^2_C(I,P')\to \Ext^2_C(S,S')$. Thus
$\Ext^2_C(I,P')\ne 0$ and consequently  $\Ext^2_C(DC,C)\ne 0$.
\end{subsection}

\begin{subsection}{A system of relations}
We wish to describe the bound quiver of a relation-extension
algebra. Let $C$ be an algebra. It is well-known that there exists a
(uniquely determined) quiver $Q_C$ and an admissible ideal $I$ of the
path algebra $kQ_C$ of $Q_C$ such that $C\isomorphe kQ_C/I$, see, for
instance, \cite{BoG}. We denote by $(Q_C)_0$ the set of points of
$Q_C$ and by $(Q_C)_1$ its set of arrows. For each point $x\in
(Q_C)_0$, we let  $e_x$ denote the corresponding primitive idempotent of
$C$, and by $S_x,P_x,I_x$ respectively,  the corresponding simple,
indecomposable projective and indecomposable injective $C$-module.

Following \cite{Bo}, we define a \emph{system of relations} for
$C\isomorphe kQ_C/I$ to be a subset $R$ of $\bigcup_{x,y \in
  (Q_C)_0}e_x I e_y$ such that $R$, but no proper subset of $R$,
generates $I$ as a two-sided ideal of $kQ_C$. Thus, for any $x,y\in
(Q_C)_0$, the elements of  $R \cap (e_x I e_y) $ are linear
combinations of paths (of length at least two) from $x$ to $y$. 
We need the following
result.
\begin{lemma}[({\cite[1.2]{Bo}})]\label{2.2} Let $C\isomorphe kQ_C/I$ be such that
  $Q_C$ has no oriented cycles and $R$ be a system of relations for $C$. Then,
  for each $x,y\in (Q_C)_0$,   the cardinality of the set $R \cap (e_x
  I e_y)$ is independent of the chosen system of relations for $C$,
  and equals   $\dim_k\Ext_C^2(S_x,S_y)$.
\end{lemma}
\end{subsection}

\begin{subsection}{The quiver of a trivial extension}\label{section 2.3}
We start with the following easy lemma.
\begin{lemma}\label{2.3}
  Let $C$ be an algebra, and $M$ be a $C$-$C$-bimodule. The quiver
  $Q_{C\ltimes M }$ of the trivial extension of $C$ by $M$ is
  constructed as follows:
\begin{enumerate}
\item $(Q_{C\ltimes M })_0=(Q_C)_0$
\item For $x,y\in (Q_C)_0$, the set of arrows in $ Q_{C\ltimes M }$
  from $x$ to $y$ equals the set of arrows in $Q_C$  from $x$ to $y$
  plus 
\[\dim_k \frac{e_x \,M \,e_y}{ e_x \,M\,(\rad C) \,e_y + e_x\, (\rad
    C)\,M \,e_y }\] 
additional arrows from $x$ to $y$.
\end{enumerate}
\end{lemma}
\begin{proof}
Since $M\subset \rad(C\ltimes M)$, the quivers of $C\ltimes M$ and of
$C$ have the same points. The arrows in the quiver of $C \ltimes M$
correspond to a $k$-basis of the vector space
 \begin{equation}\nonumber
 \rad(C \ltimes M) \; / \; \rad\!^2\,(C \ltimes M).
\end{equation}
Now, as a vector space
\begin{equation}\nonumber
 \rad(C\ltimes M) = \rad C \oplus M
\end{equation}
and since $M^2=0$ in $C \ltimes M,$
\begin{equation}\nonumber
 \rad\!^2\,(C \ltimes M) = \rad\!^2\, C \oplus \left[M\,(\rad C)+(\rad
 C)\,M\right]. 
\end{equation}
Since $\rad\!^2\, C \subset \rad C$ and $M\,(\rad C)+(\rad C)\,M \subset
 M$ and since the arrows of $Q_C$  correspond to a basis of $\rad
 C/\rad\!^2\, C$, 
 the additional arrows of  $Q_{C \ltimes M}$ correspond to a
 $k$-basis of $M / [M\,(\rad C)+(\rad C)\,M]$. 
The arrows from $x$ to $y$ are obtained upon multiplying by  $e_x$
 on the left and by $e_y$ on the right. 
\end{proof}
\end{subsection}

\begin{subsection}{The top of $\Ext^2_C(DC,C)$}\label{section 2.4}
In the situation of section \ref{section 2.3}, the $C$-$C$-bimodule
$M\,(\rad C)+(\rad C)\,M$ is the radical of $M$, and the quotient   $M
/ [M\,(\rad C)+(\rad C)\,M]$ is its top. In the case of
relation-extension algebras, we are interested in the top of
$\Ext^2_C(DC,C)$. 

\begin{lemma}\label{2.4}
  Let $C$ be an algebra of global dimension two. The top of
  the $C$-$C$-bimodule $\Ext^2_C(DC,C)$ is isomorphic to $\Ext^2_C(\soc
  DC,\tp C)$.
\end{lemma}
\begin{proof}
  The short exact sequences
\begin{equation}\nonumber
\xymatrix{0 \ar[r] & \rad C \ar[r]^{i} & C \ar[r] & \tp C \ar[r] & 0\\
 0 \ar[r] &  \soc DC \ar[r]^{j} &  DC \ar[r] &  DC/\soc DC \ar[r] &  0}
\end{equation}
where $i,j$ are the inclusions, induce a commutative diagram with
exact rows  and columns (the zeros are obtained from the condition
that the global dimension of $ C$ is two).
\begin{equation}\nonumber
\xymatrix@R=20pt@C=20pt{\acht \Ext^2_C(DC/\soc DC,\rad C)\ar[r]\ar[d]&
  {\acht \Ext^2_C(DC/\soc DC,
  C)}\ar[d]^{j^*}\ar[r] & {\acht \Ext^2_C(DC/\soc DC, \tp
  C)}\ar[d]\ar[r]&{\acht 0}\\
{\acht \Ext^2_C(DC,\rad C)} \ar[r]^{i_*}\ar[d]& {\acht\Ext^2_C(DC,C)}
  \ar[r]\ar[d]\ar@{.>}[dr]^{{  p}}&{\acht \Ext^2_C(DC,\tp
  C)}\ar[r]\ar[d] & {\acht 0} \\ 
{\acht \Ext^2_C(\soc DC,\rad C)} \ar[r]\ar[d] & {\acht \Ext^2_C(\soc DC,C)
  }\ar[r]\ar[d] & 
  {\acht \Ext^2_C(\soc DC,\tp C)} \ar[r]\ar[d] & {\acht 0}\\
{\acht 0} & {\acht 0} & {\acht 0} &}
\end{equation}
By the commutativity of the lower-right square, there exists  an epimorphism
${  p} : \Ext^2_C(DC,C) \to \Ext^2_C(\soc DC, \tp C).$
 We thus only need to show that the kernel of ${  p}$ is isomorphic to the radical
\begin{equation}\nonumber
\Ext^2_C(DC,C) \,(\rad C) + (\rad C) \,\Ext^2_C(DC,C)
\end{equation}
of the $C$-$C$-bimodule $\Ext^2_C(DC,C)$.
Now an easy diagram chasing yields
\begin{equation}\nonumber
\Ker {  p} = \im j^* + \im i_*.
\end{equation}
Thus, it suffices to prove that
\[\begin{array}{rcl}
\im i_*   =   (\rad C)\, \Ext^2_C(DC,C) & \quad \textup{and}\quad& \im j^*   =   \Ext^2_C(DC,C)\,  (\rad C) .
\end{array}\]
We only show the first equality, the second is shown similarly.
Let
\begin{equation}\nonumber
\xymatrix{0 \ar[r] & P_2 \ar[r]^{d_2} & P_1 \ar[r]^{d_1} & P_0 \ar[r]^{d_0} & DC \ar[r] & 0\\}
\end{equation}
be a projective resolution of $DC$. By definition
\begin{equation}\nonumber
\Ext^2_C(DC,C) = \Hom_C(P_2,C)/ \im \Hom_C(d_2,C) \; .
\end{equation}
We first claim that the image of the map
\begin{equation}\nonumber
i_0 = \Hom_C(P_2,i) : \Hom_C(P_2, \rad C) \to \Hom_C(P_2,C)
\end{equation}
is equal to $ (\rad C)\, \Hom_C(P_2,C).$ 
Indeed,  the product $r f $ with $r \in \rad C$ and $f \in
\Hom_C(P_2,C)$ is easily seen to factor through $\rad C$. Therfore, we have 
$ (\rad C)\, \Hom_C(P_2,C) \subset \im i_0$.
On the other hand, there is  an isomorphism  of $k$-vector spaces
\begin{equation}\nonumber
 (\rad C) \,\Hom_C (P_2,C) \isomorphe \Hom_C(P_2,\rad C). \end{equation}
Since $i_0$ is injective, this establishes our claim.

Now, the image of $i_*$ is generated by the residual classes (modulo
the image of $\Hom_C(d_2,C)$) of the products $i g $, with $g \in \Hom_C(P_2,
\rad C)$. These are the residual classes of  the elements in $\im i_0$
thus, by our claim above, the residual classes of the elements  of the
form $r f $ with $r
\in \rad C$ and $f \in \Hom_C(P_2,C)$. We deduce that $\im i_* = (\rad C)\,
\Ext^2_C(DC,C)$, as required.
 
\end{proof}

\begin{remark}\label{remark3}
The proof of this lemma  can easily be generalised to show that, for
an algebra $C$ of global dimension at most $m$, the top of the
bimodule $\Ext_C^m(DC,C)$ is equal to $\Ext_C^m(\soc DC, \tp C)$. 
\end{remark}
\end{subsection}

\begin{subsection}{The quiver of a relation-extension}\label{section 2.5}
The following theorem states that the quiver of the relation-extension
algebra is obtained from the quiver of the original algebra by adding,
for each pair of points $x,y$, one arrow from $x$ to $y$ for each
relation from $y$ to $x$. This
justifies the name ``relation-extension''.
\begin{theorem}\label{2.5}
 Let $C\isomorphe kQ_C/I$ be an algebra of global dimension
at most two, such that $Q_C$ has no oriented cycles, and let $R$ be a
system of relations for $C$. The quiver of the relation-extension
algebra $C\ltimes \Ext^2_C(DC,C)$ is constructed as follows:
\begin{itemize}
\item[(a)] $(Q_{C\ltimes \Ext^2_C(DC,C)})_0=(Q_C)_0$
\item[(b)] For $x,y\in (Q_C)_0$, the set of arrows in $Q_{C\ltimes
  \Ext^2_C(DC,C)}$ from $x$ to $y$ equals the set of arrows in $Q_C$
  from $x$ to $y$ plus $\textup{Card}\left(R\cap
  (e_y\,I\,e_x)\right)$ additional arrows.
\end{itemize}
\end{theorem}
\begin{proof}
Let $S_1, S_2, \ldots , S_n$ denote a complete set of 
representatives of the isomorphism classes of simple $C$-modules, and
set $S= \oplus_{i=1}^n S_i$. 
Since $C$ is basic, the module $S$ is isomorphic to the top
of $C$ and to the socle of $DC$. 
By Lemma \ref{2.2}, the relations of $R$ correspond to a $k$-basis of
$\Ext_C^2(S,S)$.  
By Lemma \ref{2.4}, the $C$-$C$-bimodule $\Ext_C^2(S,S)$ is isomorphic to
the top of   $\Ext_C^2(DC,C)$. 
Lemma \ref{2.3} then implies that the number of additional arrows from $x$ to
$y$ equals the $k$-dimension of the vector space $ e_x \Ext^2_C(S,S) e_y =
\Ext^2_C(S_y,S_x)$, and the result follows.
 
\end{proof}
In particular, the quiver of a non-hereditary relation-extension
algebra always contains oriented cycles.
\end{subsection}
\begin{subsection}{The indecomposable projectives}
It would be useful to know a system of relations for the
relation-extension algebra $C\ltimes\Ext^2_C(DC,C)$ starting from one
for $C$. In actual examples,
such a system is easily  obtained once we know  the
indecomposable projective modules.
In order to state the next lemma, we need a notation: for each $x\in
(Q_C)_0$, we denote by $\tilde P_x$ the corresponding indecomposable
projective $C\ltimes\Ext^2_C(DC,C)$-module. Also, we note that
$C$-modules can always be considered as $C\ltimes\Ext^2_C(DC,C)$-modules
under the standard embedding.
\begin{lemma}\label{2.6}
  Let $C$ be an algebra of global dimension at most two. Then, for
  each $x\in (Q_C)_0$, we have a short exact sequence in $\modd
  \left(C\ltimes\Ext^2_C(DC,C)\right)$ 
\begin{equation}\nonumber
0\to \Ext^2_C(DC,P_x) \to \tilde P_x \stackrel{p_x}{\to} P_x \to 0 
\end{equation}
where $p_x$ is a projective cover.
\end{lemma}
\begin{proof}
Since both $P_x$ and $\tilde P_x$ admit $S_x$ as a simple top, there
indeed exists a projective  cover morphism $p_x:\tilde P_x\to P_x$. On
the other hand, $ \Ext^2_C(DC,P_x)\isomorphe e_x \Ext^2_C(DC,C)$ is
clearly a submodule of the $C\ltimes  \Ext^2_C(DC,C)$-module $\tilde
P_x$. The result then follows from the isomorphism of $k$-vector
spaces 
\begin{equation}\nonumber
\tilde P_x = e_x\left(C\ltimes \Ext^2_C(DC,C)\right)
\isomorphe P_x \oplus  \Ext^2_C(DC,P_x).
\end{equation}
 
\end{proof} 
\end{subsection} 
\begin{subsection}{An example}
\begin{example}\label{example 2.7}
Let $C$ be given by the quiver 
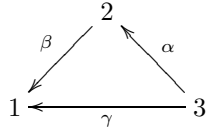
\[
\xymatrix{&2\ar[dl]_\zb\\ 1&&3\ar[ll]^\zg \ar[ul]_\za}
\]
bound by the relation $\za\zb=0$. Thus
\begin{equation}\nonumber
C_C= 
1\,\oplus\begin{array}{c} \xymatrix@R=0pt@C=0pt{2\\1}\end{array} \oplus
\begin{array}{c}   \xymatrix@R=0pt@C=0pt{&3\\ 1&& 2} \end{array}
\quad \textup{  and  } \quad
(DC)_C= \begin{array}{c}\xymatrix@R=0pt@C=0pt{2&& 3\\& 1}\end{array}\oplus
\begin{array}{c}\xymatrix@R=0pt@C=0pt{ 3\\ 2}\end{array} \oplus \, 3
\end{equation}
where the indecomposable projectives and injectives are represented by
their Loewy series. It is easily seen that the global dimension of $C$
is two. By Theorem \ref{2.5}, the quiver of   $C\ltimes\Ext^2_C(DC,C)$
is obtained by adding to $Q_C$ a single arrow $\zd:1\to 3$.
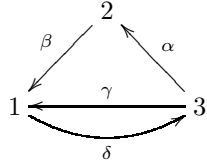
\[\xymatrix{&2\ar[dl]_\zb\\ 1\ar@/_12pt/[rr]_\zd&&3\ar[ll]_\zg
  \ar[ul]_\za}\]
We now compute the new indecomposable projective modules. A simple
calculation yields 
\[\begin{array}{rclcrcl}
  \Ext^2_C(I_3,P_1)&\isomorphe& k &\quad,\quad&
  \Ext^2_C(I_3,P_3)&\isomorphe& k\\ \\
  \Ext^2_C(I_1,P_1)&\isomorphe& k &,&   \Ext^2_C(I_1,P_3)&\isomorphe
  &k. 
\end{array}\]
Since the projective dimension of $I_2$ is one and the injective
dimension of $P_2$ is also one, this yields
$\dim_k\Ext_C^2(DC,C)=4$. 
Using Lemma \ref{2.6},
 we get the new indecomposable projectives
\[\begin{array}{c}\xymatrix@R=0pt@C=0pt{1\\3\\1}\end{array}\qquad,\qquad
 \begin{array}{c}\xymatrix@R=0pt@C=0pt{&3\\  1&&2\\3\\1}\end{array}\qquad,\qquad 
\begin{array}{c}\xymatrix@R=0pt@C=0pt{2\\1}\end{array}\ .\]

Thus, a system of relations for the relation-extension algebra is
$\za\zb=0,\zd\za=0,\zb\zd=0$ and $\zd\zg\zd=0$.
\end{example}
\end{subsection}
\end{section}


\begin{section}{Cluster-tilted algebras}
\begin{subsection}{Preliminaries}
Let $A$ be a hereditary algebra. The cluster category $\CA$ of $A$ is
defined as follows. Let $F$ denote the automorphism of $\DA$ defined
as the composition $\tau^{-1}_{\DA}\, [1]$, where $\tau^{-1}_{\DA}$ is
the Auslander-Reiten translation in $\DA$, and $[1]$ is the shift
functor. Then $\CA$ is the quotient category $\DA/F$. Its objects are
the $F$-orbits $\tilde X=(F^iX)_{i\in \mathbb{Z}}$, where $X$ is an object
in $\DA$. The set of morphisms from  $\tilde X=(F^iX)_{i\in \mathbb{Z}}$ to
$\tilde Y=(F^iY)_{i\in \mathbb{Z}}$ in $\CA$ is given by
\begin{equation}\nonumber
\Hom_{\CA}(\tilde X,\tilde Y) = \bigoplus_{i\in \mathbb{Z}} \Hom_{\DA}(X,F^iY).
\end{equation} 
It is shown in \cite{K}, that $\CA$ is a triangulated
category. Furthermore, the canonical functor $\DA\to \CA$ is a functor
of triangulated categories. We refer to \cite{BMRRT} for facts about
the cluster category.

An object $\tilde T$ in $\CA$ is called a \emph{tilting object}
provided $\Ext^1_{\CA}(\tilde T,\tilde T)=0$ and the number of
isomorphism classes of indecomposable summands of $\tilde T$ equals
the number of isomorphism classes of simple $A$-modules (that is, the
number of points in the quiver of $A$). The algebra of endomorphisms
$\tilde C=\End_{\CA}(\tilde T)$ is then called a \emph{cluster-tilted
  algebra} \cite{BMR1}. 

Cluster-tilted algebras  may also be expressed in terms of modules. We
recall that an $A$-module $T$ is called a tilting module provided
$\Ext^1_A(T,T)=0$ and the number of
isomorphism classes of indecomposable summands of $ T$ equals
the number of isomorphism classes of simple $A$-modules. Denoting by
$\tilde T $ the $F$-orbit of $T$, we have 
the following theorem.
\begin{theorem}[({\cite[3.3]{BMRRT}})]\label{thm bmrrt}
Let $\tilde C$ be a cluster-tilted algebra, then there exist a
hereditary algebra $A$ and a tilting $A$-module $T$ such that $\tilde
C\isomorphe \End_{\CA}(\tilde T)$.
\end{theorem} 

We further recall that the endomorphism algebra of a tilting module
over a hereditary algebra is called a \emph{tilted algebra}, see, for
instance, \cite{R}. We need the following result.
\begin{theorem}[({\cite{H}})]\label{happelthm} Let $A$ be a hereditary
  algebra, $T$ be a tilting $A$-module  
 and $C=\End_A(T)$ be the corresponding
  tilted algebra. Then
\begin{itemize}
\item[(a)] The derived functor
$\textup{RHom}_A(T,-): \DA\to \DC$ is an equivalence of
  categories which maps the $A$-module $T$ to the $C$-module $C$. 
\item[(b)] $\textup{RHom}_A(T,-)$ commutes with the Auslander-Reiten
  translations and the shifts in the respective categories.
\end{itemize}
\end{theorem}
\end{subsection} 
\begin{subsection}{Cluster-tilted algebras are trivial extensions} For any object $X$ in $\DA$, the $k$-vector space
  $\Hom_{\DA}(X,FX)$ has a natural structure of
  $\End_{\DA}(X)$-$\End_{\DA}(X)$-bimodule under the action 
\[\begin{array}{rcl} 
\End_{}(X)\times\Hom_{}(X,FX)\times\End_{}(X) &\to&
\Hom_{}(X,FX)\\
(u,f,v)&\mapsto &Fu\circ f \circ v
\end{array} \]
The following lemma is proved in \cite[3.1]{B}. We include a simple
proof for the convenience of the reader.

\begin{lemma}\label{3.2} Let $\tilde C$ be a cluster tilted
  algebra. Then, for each hereditary algebra $A$ and tilting
 $A$-module $T$  such that $\tilde C = \End_{\CA}(\tilde T)$, we have
\begin{equation} \nonumber
\tilde C \isomorphe \End_A(T)\ltimes
  \Hom_{\DA}(T, FT).
\end{equation}
\end{lemma}
\begin{proof} By definition of  $\CA$, we have  
\begin{equation} \tilde C = \End_{\CA}(\tilde T) =
  \oplus_{i\in\mathbb{Z}}\ \Hom_{\DA}(T,F^i\,T )
\nonumber
\end{equation}
as $k$-vector spaces, and the multiplication is given by 
\begin{equation}\nonumber
(g_i)_{i\in\mathbb{Z}}(f_j)_{j\in \mathbb{Z}}= \left(\sum_{i+j=l} F^j\,g_i\circ
  f_j\right)_{l\in\mathbb{Z}}.
\end{equation} 
Since $A$ is hereditary, then, for any two $A$-modules $M$ and $N$,
we have that $\Hom_{\DA}(M,N[i])=0$ for all $i\ge 2$. Therefore, as a
$k$-vector space  
   \begin{equation}\nonumber
\tilde C=\End_{\CA}(\tilde T) \  =\ 
  \Hom_{\DA}(T,T)\oplus \Hom_{\DA}(T,FT).
\end{equation}
The multiplication of two elements $f,g \in \End_{\CA}(\tilde T) $ is  given as
follows.
Assume $f=(f_0\,,f_1)$ and $ g=(g_0\,,g_1)$, with
$f_0,g_0\in\Hom_{\DA}(T,T)$ and $f_1,g_1 \in 
\Hom_{\DA}(T,FT)$, then, since  $Fg_1\circ f_1=0$,
\begin{equation}\nonumber
 gf =(g_0\circ f_0\ ,\, Fg_0\circ f_1+f_0 \circ g_1).
\end{equation}  
 In view of the bimodule structure of
$\Hom_{\DA}(T,FT)$ defined above, this shows indeed that  $\tilde
C=\End_{\CA}(\tilde T)$ 
is the trivial extension of $ \End_{\DA}(T)=\End_A(T)$ by the bimodule
 $\Hom_{\DA}(T,FT).$ 
 
\end{proof}

Since the algebra $\End_A(T)$ of the lemma is tilted, any
  cluster-tilted algebra  is a trivial extension of a tilted algebra.
However, the hereditary algebra $A$ and the $A$-module $T$ above are not
  unique. Therefore, one cannot apply directly the lemma to construct
  a map from cluster tilted algebras to
  tilted algebras.

\end{subsection}


\begin{subsection}{The main result}
We are now able to prove the main theorem of this section.

\begin{theorem}\label{3.3}
An algebra $\tilde C$ is cluster-tilted if and
    only if there exists a tilted algebra $C$ such that $\tilde C$ is
    the relation-extension of $C$.
\end{theorem}
\begin{proof} 
Let $C$ be a tilted algebra. Then there exist a hereditary algebra $A$
  and a tilting $A$-module 
  $T$ such that $C=\End_A(T)$. Let $\tilde T$  denote as usual  the
  $F$-orbit of $T$ in $\DA$. Then $\tilde C = \End_{\CA}(\tilde T)$
  is a 
  cluster-tilted algebra.
By Lemma \ref{3.2}, we have
\begin{equation}
  \tilde C = \End_{\DA}(T)\ltimes
  \Hom_{\DA}(T,FT).
\end{equation}
By Theorem \ref{happelthm}, the derived functor $\textup{RHom}_A(T,-)$
induces $C$-$C$-bimodule isomorphisms
\begin{equation}\nonumber
\End_{\DA}(T)\isomorphe\End_{\DC}(C)\isomorphe\End_C(C)\isomorphe C
\end{equation}
and \begin{equation}\nonumber
 \Hom_{\DA}(T,FT) \isomorphe  \Hom_{\DC}(C,F'C)
\end{equation}
where $F'=\tau^{-1}_{\DC}\,[1]$ is the functor corresponding to $F$ in
the derived category $\DC$. Thus we get
\begin{equation}\nonumber
 \tilde C \isomorphe C \ltimes  \Hom_{\DC}(C,F'C).
\end{equation} 
Moreover, we have the following sequence of $C$-$C$-bimodule isomorphisms
\[\begin{array}{rcl}
\Hom_{\DC}(C,F'C)
 &\ \isomorphe \ & \Hom_{\DC}(\,\tau_{\DC} \,C \,[1] \,,\, C\,
 [2]\,)\\ 
 &\isomorphe & \Hom_{\DC}(\,DC\,, \,C [2]\,)\\
 &\isomorphe &\Ext^2_C(DC,C),
\end{array}\]
where the first  is obtained by applying to both arguments the
automorphism $\tau_{\DC}\,[1]$, the second uses the fact that
$\tau_{\DC} C \isomorphe DC[-1]$ and 
the third is a property of the derived category.
This shows that the relation-extension $C\ltimes \Ext^2_C(DC,C)$ is a
cluster-tilted algebra. Finally, by Lemma \ref{3.2}, every
cluster-tilted algebra  is obtained in this way.
\end{proof}
\end{subsection}

\begin{subsection}{Remarks and examples} 
\begin{itemize}
\item[(a)] Since the quiver of a tilted algebra has no oriented
  cycles, it follows directly from Theorem \ref{3.3} and Theorem
  \ref{2.5} that we have a construction for the quiver of a
  cluster-tilted algebra  $\tilde C$ starting from the quiver of a
  tilted algebra $C$. This construction is easily seen to generalise
  the one in \cite[4.1]{BRS} and, thus, can be used to relate the
  Happel-Vossieck list of tame concealed algebras \cite{HV} with
  Seven's list of minimal infinite cluster quivers \cite{S}.
\item[(b)] A different description, inspired from \cite{HW}, of the
  relation-extension algebra is sometimes useful. Consider the
  following doubly infinite matrix algebra 
\begin{equation}\nonumber
\hat C  \quad = \quad \left[
  \begin{array}{cccccccccc}
  \ddots &&&\ 0\ \\
&\ C_{i-1}\ \\
&M_i&\ C_i\ \\
&&M_{i+1}&\ C_{i+1}\ \\
&\ 0\ &&&\ddots
  \end{array}
 \right] 
\end{equation}
where matrices are assumed to have only finitely many non-zero
coefficients, $C_i=C$ and $M_i=\Ext_C^2(DC,C)$ for all $i\in \mathbb{Z}$,
all the remaining coefficients are zero. The addition is the usual
addition of matrices while the multiplication is induced from the
bimodule structure of $\Ext_C^2(DC,C)$ and the zero map
$\Ext_C^2(DC,C)\otimes_C\Ext_C^2(DC,C)\to 0$. Clearly, $\hat C$ is a
Galois covering of $C\ltimes \Ext_C^2(DC,C)$ with group $\mathbb{Z}$: the
identity maps $C_i\to C_{i+1}$, $M_i\to M_{i+1}$ induce an
automorphism $\eta$ of $\hat C$ and $\hat C / \eta \isomorphe C
\ltimes \Ext_C^2(DC,C)$.
\item[(c)] As observed before, different tilted algebras $C$ may
  correspond to the same cluster-tilted algebra  $\tilde C$ (thus, the
  surjective  map $C\mapsto \tilde C$ is not injective). We give an example
  of such an occurrence.
\begin{example}
Let $C_1$ be given by the quiver
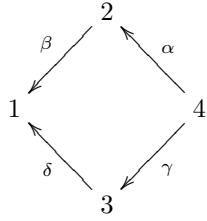
\begin{equation}\nonumber
\xymatrix{&2\ar[ld]_\zb \\ 1&&4\ar[lu]_\za\ar[ld]^\zg\\ &3\ar[lu]^\zd}
\end{equation}
bound by $\za\zb=\zg\zd$. This is a tilted algebra of Dynkin type
$D_4$, and the corresponding cluster-tilted (relation-extension)
algebra  $\tilde C_1$ is given by the quiver
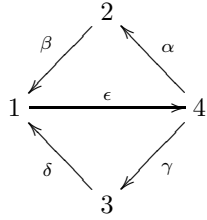
\begin{equation}\nonumber
\xymatrix{&2\ar[ld]_\zb \\ 1\ar[rr]^\ze&&4\ar[lu]_\za\ar[ld]^\zg\\
  &3\ar[lu]^\zd} 
\end{equation}
bound by $\za\zb=\zg\zd$, $\zb\ze=0$, $\zd\ze=0$, $\ze\za=0$,
$\ze\zg=0$. 
Let now $C_2$ be the tilted algebra given by the quiver 
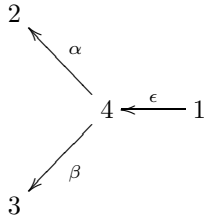
\begin{equation}\nonumber
\xymatrix{2\\ &4\ar[ld]^\zb\ar[lu]_\za&1\ar[l]_\ze\\ 3}
\end{equation}
bound by $\ze\za=0$, $\ze\zb=0$. Then it is easily seen that $\tilde
C_1=\tilde C_2$.
\end{example}
\item[(d)] Not surprisingly, it is possible that  $C$ is
 representation-finite whereas $\tilde C$ is  representation-infinite: it
  suffices to have two points $x,y\in (Q_C)_0$ such that
  $\dim_k\Ext_C^2(I_y,P_x)>1.$ We give an example of such a situation.
\begin{example}
Let $C$ be given by the quiver
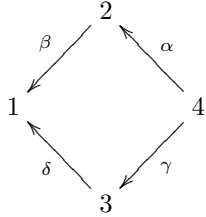
\begin{equation}\nonumber
\xymatrix{&2\ar[ld]_\zb \\ 1&&4\ar[lu]_\za\ar[ld]^\zg\\ &3\ar[lu]^\zd}
\end{equation}
bound by $\za\zb=0$, $\zg\zd=0$. This is a representation-finite
tilted algebra of euclidean type $\tilde A_3$. The injective
resolution 
\begin{equation}\nonumber 
  0\to P_1\to I_1\to I_2\oplus I_3\to I_4\oplus I_4\to 0
\end{equation} 
shows that   $\dim_k\Ext_C^2(I_4,P_1)=2.$ 
The corresponding cluster-tilted 
algebra  $\tilde C$ is given by the quiver
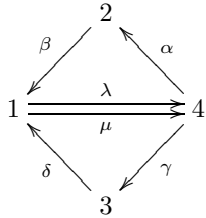
\begin{equation}\nonumber 
\xymatrix{&2\ar[ld]_\zb \\
  1\ar@<2pt>[rr]^\zl\ar@<-2pt>[rr]_\mu&&4\ar[lu]_\za\ar[ld]^\zg\\ 
  &3\ar[lu]^\zd} 
\end{equation} 
bound by $\za\zb=0$, $\zg\zd=0$,  $\zd\zl=0$, $\zl\zg=0$,
$\zb\mu=0$, $\mu\za=0$. 
The indecomposable projective $\tilde C$-modules are given by
\begin{equation}\nonumber
\begin{array}{c}\xymatrix@R=0pt@C=0pt{&1\\4&&4\\2&&3}\end{array}
\qquad,\qquad
\begin{array}{c}\xymatrix@R=0pt@C=0pt{2\\1\\4\\2}\end{array}
\qquad,\qquad
\begin{array}{c}\xymatrix@R=0pt@C=0pt{3\\1\\4\\3}\end{array}
\qquad,\qquad
\begin{array}{c} \xymatrix@R=0pt@C=0pt{&4\\2&&3}\end{array}
\end{equation} 
Clearly, $\tilde C$ is representation-infinite. 
\end{example}
\item[(e)] The relation-extension algebra in Example \ref{example 2.7}
  is not a cluster-tilted algebra. This follows from the fact that
  cluster-tilted algebras contain no oriented cycles of length two.

\end{itemize}
  
\end{subsection} 

\end{section}

\vfill

\noindent  I. Assem\\
   D\'epartement de Math\'ematiques,\\
  Universit\'e de Sherbrooke,\\
  Sherbrooke (Qu\'ebec), J1K 2R1, Canada\\
  {ibrahim.assem@usherbrooke.ca}

\bigskip

\noindent  T. Br\"ustle\\
 D\'epartement de Math\'ematiques,\\
  Universit\'e de Sherbrooke, \\
Sherbrooke (Qu\'ebec), J1K 2R1  Canada\\
   {thomas.brustle@usherbrooke.ca}\\
and\\ Department of Mathematics,\\
 Bishop's University,\\
  Lennoxville, (Qu\'ebec),
  J1M 1Z7,
 Canada\\
      {tbruestl@ubishops.ca}

\bigskip

\noindent R. Schiffler\\
Department of Mathematics and  Statistics,\\
 University of Massachusetts at Amherst,\\
 Amherst, MA
  01003-9305,
 USA \\
{schiffler@math.umass.edu}

\end{document}